# Several Fixed Point Theorems
# on Partially Ordered Banach Spaces and Applications

**Jinlu Li**
Department of Mathematics
Shawnee State University
Portsmouth, Ohio 45662
USA

**Abstract**

In this paper, we prove several fixed point theorems on both of normal partially ordered Banach spaces and regular partially ordered Banach spaces by using the normality, regularity, full regularity, and chain -complete property. Then, by applying these theorems, we provide some existence and uniqueness of solutions to some integral equations. We also prove the solvability of some equilibrium problems in Banach spaces.

**Keywords:** universally inductive poset; partially ordered Banach space; fixed point

**MSC:** 06F30; 06F30; 45G10; 91A06; 91A18

## 1. Introduction

In the traditional fixed point theory, the underlying spaces are topological spaces and the considered mappings must satisfy a certain type of continuity to insure the existence of fixed point. Tarski's Fixed Point Theorem on chain-complete lattice for single-valued mappings (see [12]) initiated a new custom in fixed point theory, in which there are some ordering relations on the underlying spaces, such as, preorder, partial order, or lattice, and the underlying spaces are not required to be equipped with topological structure. To guarantee the existence of fixed point, the considered mappings should satisfy some order-monotonic conditions and it is unnecessary for them to have any continuity property.

Based on the fact that Banach spaces are the fundamental underlying spaces on linear and nonlinear analysis, it leads us to consider the following problem: if a Banach space is equipped with an ordering structure, partial order or lattice, this Banach space becomes a partially ordered Banach space. Then, when we solve some problems on this Banach space, in addition to the topological structure and the algebraic structure, the ordering structure will provide a new powerful tool. This important idea has been widely used in solving integral equations ([3-5], [7-8], [10], [11], [13]), vector variational inequalities ([7]), nonlinear fractional evolution equations ([14]), Nesh equilibrium problems ([2], [6], [15]), etc.

For example, by applying Theorems 3.10 in [9], we can obtain the following result.

**Theorem**. *Let $(X, \|\cdot\|, \succcurlyeq)$ be a partially ordered reflexive Banach space and let D be a bounded closed convex subset of X. Let $F: D \to 2^D\setminus\{\emptyset\}$ be a set-valued mapping satisfying the following three conditions*:

A1. *F is $\succcurlyeq$-increasing upward*;
A2. *F(x) is a closed and convex subset of D, for every $x \in D$*;
A3. *There are elements $y \in D$ and $v \in F(y)$ with $y \preccurlyeq v$*.

*Then*

(i) *$(\mathcal{F}(F), \succcurlyeq)$ is a nonempty inductive poset*.

(ii) *$(\mathcal{F}(F) \cap [y), \succcurlyeq)$ is a nonempty inductive poset*.

In this theorem, the underlying space is a partially ordered Banach space. The set-valued mapping $F$ satisfies some order-increasing upward conditions A1-A3. Of cause, this mapping $F$ does not have any continuity property. Then the existence of fixed point of $F$ is guaranteed. Furthermore, it provides the inductive properties of the set of fixed points of $F$. No doubt, Tarski's Fixed Point Theorem, the theorem listed above, and other fixed point theorems on posets (see [1], [6], [9], [12]) provide useful tools in analysis in ordered sets, such as the equilibrium problems with incomplete preferences (see [1], [6], [15]).

To solve more difficult problems in analysis on partially ordered Banach spaces, it is important to develop more fixed point theorems by using both of the properties of the topologies and the properties of the partial orders. This is the goal of this paper.

In this paper, we introduce the concept of $\delta$-distance between bounded sets (see section 3) and condition $(H^1)$ (see section 5) for mappings in Banach spaces. Then by these concepts and based on the results obtained in [5] and [8] about the connections between the normality, regularity and chain-complete properties, we prove several fixed point theorems for both single-valued and set-valued mappings on partially ordered Banach spaces.

## 2. Preliminaries

Let $(X, \|\cdot\|)$ be a Banach space and $K$ a nonempty closed convex cone of $X$. It is well known that a partial order $\succcurlyeq$ on $X$ can be induced by $K$ as follows.

$$x \succcurlyeq y \text{ if and only if } x - y \in K, \text{ for all } x, y \in X.$$

Then $(X, \|\cdot\|, \succcurlyeq)$ is a partially ordered Banach space which satisfies that, for any $u, w \in X$, the following $\succcurlyeq$-intervals are closed with respect to the $\|\cdot\|$-topology:

$$[u) = \{x \in X : x \succcurlyeq u\} \text{ and } (w] = \{x \in X : x \preccurlyeq w\}.$$

Throughout this paper, unless otherwise is stated, we say that a partially ordered Banach space $(X, \|\cdot\|, \succcurlyeq)$ is chain-complete if every $\succcurlyeq$-upper bounded chain in $X$ has the smallest $\succcurlyeq$-upper bound in $X$; we say that $(X, \|\cdot\|, \succcurlyeq)$ is bi-chain-complete if every $\succcurlyeq$-bounded chain in $X$ has both smallest $\succcurlyeq$-upper bound and greatest $\succcurlyeq$-lower bound.

If there is a constant $\lambda > 0$ such that

$$0 \preccurlyeq x \preccurlyeq y \text{ implies that } \|x\| \leq \lambda \|y\|,$$

then $\succcurlyeq$ is said to be normal and $(X, \|\cdot\|, \succcurlyeq)$ is called a normal partially ordered Banach space. The minimum value of $\lambda$ satisfying the above inequality is called the normal constant of $\succcurlyeq$.

If every $\succcurlyeq$-upper bounded and $\succcurlyeq$-increasing sequence $\{x_n\}$ of $X$ is an $\|\cdot\|$-convergent sequence, then $\succcurlyeq$ is said to be regular and $(X, \|\cdot\|, \succcurlyeq)$ is called a regular partially ordered Banach space. If every $\|\cdot\|$-bounded and $\succcurlyeq$-increasing sequence $\{x_n\}$ of $X$ is an $\|\cdot\|$-convergent sequence, then $\succcurlyeq$ is said to be fully regular and $(X, \|\cdot\|, \succcurlyeq)$ is called a fully regular partially ordered Banach space.

We recall some results below from [8-10] for easy reference.

**Theorem 2.2.2** in [10]. *Let $(X, \|\cdot\|, \succcurlyeq)$ be a partially ordered Banach space. Then*

$$\succcurlyeq \text{ is fully regular} \Longrightarrow \succcurlyeq \text{ is regular} \Longrightarrow \succcurlyeq \text{ is normal}.$$

**Lemma 2.3** in [8]. *Let $\{x_n\}$ be an $\succcurlyeq$-increasing sequence (a chain) in a partially ordered topological space $(X, \tau, \succcurlyeq)$. If $x_n \to x$, as $n \to \infty$, then*

$$\vee \{x_n\} = x.$$

According to our definition of chain-completeness and bi-chain-completeness in this paper, from Theorem 3.6 in [8], we have

**Proposition 2.1**. *Every regular partially ordered Banach space is bi-chain-complete.*

**Corollary 3.7** in [11]. *Let $(X, \|\cdot\|, \succcurlyeq)$ be a regular partially ordered Banach space. Let $D$ be a closed inductive subset of $X$. Let $F: D \to D$ be an $\succcurlyeq$-increasing single-valued mapping. Suppose that there is $x_0 \in D$ satisfying $x_0 \preccurlyeq Fx_0$. Then*

  (a) $\mathcal{F}(F)$ *is a nonempty inductive subset of $D$;*
  (b) $\mathcal{F}(F) \cap [x_0)$ *is a nonempty inductive subset of $D$.*

## 3. $\delta$-distance between bounded sets and fixed point theorems of set-valued mappings on normal partially ordered Banach spaces

### 3.1. $\delta$-distance and its properties

**Definition 3.1.** Let $A$, $B$ be nonempty bounded subsets of a Banach space $(X, \|\cdot\|)$. The $\delta$-distance between $A$ and $B$ is defined by

$$\delta(A, B) = \max\{\sup\{\inf\{\|x-y\|: y \in B\}: x \in A\}, \sup\{\inf\{\|x-y\|: x \in A\}: y \in B\}\}.$$

**Lemma 3.2**. $\delta$ has the following properties:

  (i) $\delta(A, A) = 0$;
  (ii) $\delta(\{x\},\{y\}) = \|x-y\|$;
  (iii) $\delta(A, B) = \delta(B, A)$;

(iv) $\delta(A', B') = \delta(A, B)$, where $A'$ is the $\|\cdot\|$-closure of $A$;

(v) $\delta(A, B) = 0$, if and only if $A' = B'$.

*Proof.* Parts (i), (ii) and (iii) are clear to see.

Proof of part (iv). For any $a \in A'$ and $b \in B'$, there is $\{y_n\} \subseteq B$ with $\|y_n - b\| \to 0$, as $n \to \infty$. It implies

$$\inf\{\|a-y\|: y \in \{y_n\} \cup \{b\}\} = \inf\{\|a-y\|: y \in \{y_n\}\}. \tag{1}$$

From (1), we obtain

$$\inf\{\|a-y\|: y \in B'\} = \inf\{\|a-y\|: y \in B\}, \text{ for any } a \in A'. \tag{2}$$

On the other hand, for any given $\epsilon > 0$ and, for arbitrary $a \in A'$, there is $x \in A$ such that

$$\|a-x\| < \epsilon. \tag{3}$$

By (2) and (3), we have

$$\begin{aligned}&\inf\{\|a-y\|: y \in B'\} \\ &= \inf\{\|a-y\|: y \in B\} \\ &\leq \inf\{\|x-y\| + \|a-x\|: y \in B\} \\ &< \inf\{\|x-y\|: y \in B\} + \epsilon.\end{aligned}$$

It implies

$$\begin{aligned}&\sup\{\inf\{\|a-y\|: y \in B'\}: a \in A'\} \\ &\leq \sup\{\inf\{\|x-y\|: y \in B\}: x \in A\} + \epsilon.\end{aligned} \tag{4}$$

On the other hand,

$$\begin{aligned}&\sup\{\inf\{\|x-y\|: y \in B'\}: x \in A'\} \\ &= \sup\{\inf\{\|x-y\|: y \in B\}: x \in A'\} \\ &\geq \sup\{\inf\{\|x-y\|: y \in B\}: x \in A\}.\end{aligned} \tag{5}$$

By (4) and (5), we obtain

$$\sup\{\inf\{\|x-y\|: y \in B'\}: x \in A'\} = \sup\{\inf\{\|x-y\|: y \in B\}: x \in A\}. \tag{6}$$

We can similarly show that

$$\sup\{\inf\{\|x-y\|: x \in A'\}: y \in B'\} = \sup\{\inf\{\|x-y\|: x \in A\}: y \in B\}. \tag{7}$$

Combining (6) and (7) implies

$$\delta(A', B') = \delta(A, B).$$

*Proof of part* (v). "$\Longrightarrow$" Suppose that $\delta(A, B) = 0$. Assume, by the way of contradiction, that $A' \neq B'$, and, without loss of generality, there is $a \in A' \setminus B'$. So $a \in X \setminus B'$, which is open. It implies that there is $2r > 0$ such that

$$\{z \in X: \|z-a\| \leq 2r\} \subseteq X \setminus B'.$$

Since $a \in A'$, there is $c \in A$ with $\|c-a\| \le r$, it implies

$$\inf\{\|c-y\|: y \in B\} \ge \inf\{\|a-y\| - \|c-a\|: y \in B\} \ge 2r - r = r.$$

It follows that

$$\delta(A, B) \ge \sup\{\inf\{\|x-y\|: y \in B\}: x \in A\} \ge \inf\{\|c-y\|: y \in B\} \ge r > 0.$$

It contradicts to the assumption that $\delta(A, B) = 0$.

The part "$\Longleftarrow$" part (v) follows from part (iv) and (i) immediately. □

**Definition 3.3.** Let $(X, \|\cdot\|)$ be a Banach space and $D$ a nonempty closed subset of $X$. A set-valued mapping $T: D \to 2^X \setminus \{\emptyset\}$ with bounded values is said to be $\delta$-continuous, whenever for any Cauchy sequence $\{x_n\}$,

$$x_n \to x, \text{ as } n \to \infty \text{ and } x \in D \implies \delta(Tx_n, Tx) \to 0, \text{ as } n \to \infty.$$

$T$ is said to be $\delta$-compact if

(i) it is $\delta$-continuous;
(ii) for any bounded subset $C \subseteq D$, the set $\cup_{x \in C} Tx$ is relatively compact.

**Observation 3.4.** From Property (ii) of $\delta$-distance, it is clear to see that a single-valued mapping, as a special case of set-valued mappings with singleton values, is continuous (compact) if and only if it is $\delta$-continuous ($\delta$-compact).

**Lemma 3.5.** *Let $(X, \|\cdot\|)$ be a Banach space and $D$ a closed subset of $X$. Let $T: D \to 2^X \setminus \{\emptyset\}$ be a $\delta$-continuous mapping with closed and bounded values. Let $\{x_n\} \subseteq D$ be a Cauchy sequence with limit $x \in D$. Suppose that one of the following conditions holds*:

(i) $x_{n+1} \in Tx_n$, for $n = 1, 2, \ldots$ ;
(ii) $x_n \in Tx_n$, for $n = 1, 2, \ldots$ .

*Then $x \in Tx$.*

*Proof.* Suppose that condition (i) holds. Since $T: D \to 2^X \setminus \{\emptyset\}$ is a $\delta$-continuous mapping and $x_n \to x$, as $n \to \infty$, then

$$\delta(Tx_n, Tx) \to 0, \text{ as } n \to \infty. \qquad (8)$$

Assume, by the way of contradiction, that $x \notin Tx$. Since $Tx$ is closed, then there is $3r > 0$ such that

$$\{z \in X: \|z-x\| \le 3r\} \subseteq X \setminus Tx.$$

It implies that

$$\|z-x\| \ge 3r, \text{ for all } z \in Tx. \qquad (9)$$

From (8) and $x_n \to x$, as $n \to \infty$, for the positive number $r$, there is a large $N$, such that

$$\|x_{n+1}-x\| < r, \text{ for all } n \ge N, \qquad (10)$$

and

$$\delta(Tx_n, Tx) < r, \text{ for all } n \geq N. \tag{11}$$

Since $x_{n+1} \in Tx_n$, then, for every $n \geq N$, from (11), we have

$$\inf\{\|x_{n+1}-z\|: z \in Tx\} \leq \sup\{\inf\{\|y-z\|: z \in Tx\}: y \in Tx_n\} \leq \delta(Tx_n, Tx) < r. \tag{12}$$

From (12), it implies that, for any fixed $n \geq N$, there is $b_n \in Tx$ such that

$$\|x_{n+1}-b_n\| < r, \tag{13}$$

From (10), (9) and (13), for all $n \geq N$, we get

$$r > \|x_{n+1}-x\| \geq \|x-b_n\| - \|x_{n+1}-b_n\| > 3r - r = 2r.$$

It is a contradiction. We can similarly prove the case if condition (ii) holds. □

### 3.2. Some fixed point theorems for set-valued mappings on normal partially ordered Banach spaces

We recall some order-monotone concepts below which are used in the sequel of this section and the following sections.

Let $(X, \succcurlyeq)$, $(U, \succcurlyeq^U)$ be posets and $T: X \to 2^U\setminus\{\emptyset\}$ a set-valued mapping. $T$ is said to be isotone, or $\succcurlyeq$-increasing upward, if, $x \preccurlyeq y$ in $X$ implies, for any $z \in Tx$, there is a $w \in Ty$ such that $z \preccurlyeq^U w$. $T$ is said to be $\succcurlyeq$-increasing downward, if $x \preccurlyeq y$ in $X$ implies, for any $w \in Ty$, there is a $z \in Tx$ such that $z \preccurlyeq^U w$. If $T$ is both $\succcurlyeq$-increasing upward and $\succcurlyeq$-increasing downward, then $T$ is said to be $\succcurlyeq$-increasing.

In particular, a single-valued mapping $F$ from a poset $(X, \succcurlyeq)$ to a poset $(U, \succcurlyeq^U)$ is said to be $\succcurlyeq$-increasing whenever, for $x, y \in X$, $x \preccurlyeq y$ implies $F(x) \preccurlyeq^U F(y)$. An $\succcurlyeq$-increasing mapping $F: X \to U$ is said to be strictly $\succcurlyeq$-increasing whenever $x \prec y$ implies $F(x) \prec^U F(y)$.

Let $(X, \|\cdot\|, \succcurlyeq)$ be a partially ordered Banach space. Let $D$ be a subset of $X$. Let $T: D \to 2^X\setminus\{\emptyset\}$ be a set-valued mapping. For a point $x \in D$, if $x \in Tx$, then $x$ is called a fixed point of $T$. The collection of all fixed points of $F$ is denoted by $\mathcal{F}(T)$.

A nonempty subset $A$ of a poset $(X, \succcurlyeq)$ is said to be universally inductive in $X$ whenever, for any given chain $\{x_\alpha\} \subseteq X$, if every element $x_\beta \in \{x_\alpha\}$ has an $\succcurlyeq$-upper cover in $A$, then $\{x_\alpha\}$ has an $\succcurlyeq$-upper bound in $A$. Some useful universally inductive subsets in posets are provided in [9], which are listed as lemmas below for easy reference.

**Lemma 3.7 [9].** *Every inductive subset $A$ in a chain complete poset such that $A$ has a finite number of maximal elements is universally inductive.*

**Lemma 3.8 [9].** *Every nonempty compact subset of a partially ordered Hausdorff topological space is universally inductive.*

**Theorem 3.6.** *Let $(X, \|\cdot\|, \succcurlyeq)$ be a normal partially ordered Banach space and $D$ a closed inductive subset of $X$. Let $T: D \to 2^D \setminus \{\emptyset\}$ be an isotone and $\delta$-compact mapping with closed and bounded values. Suppose that there are points $x_0 \in D$, $x_1 \in Tx_0$ satisfying $x_0 \preccurlyeq x_1$. Then*

(a) *$\mathcal{F}(T)$ is a nonempty chain-complete subset of $D$;*

(b) *if $D$ is bi-inductive, then $\mathcal{F}(T)$ is a nonempty bi-chain-complete subset of $D$.*

*Moreover,*

(a′) *$T$ has an $\succcurlyeq$-maximal fixed point;*

(b′) *if $D$ is bi-inductive, then $T$ has both $\succcurlyeq$-maximal and $\succcurlyeq$-minimal fixed points.*

*Proof.* We first prove (a). For the given points $x_0 \in D$, $x_1 \in Tx_0$, by the condition $x_0 \preccurlyeq x_1$ and from the isotone property of $T$, there is $x_2 \in Tx_1$ satisfying

$$x_0 \preccurlyeq x_1 \preccurlyeq x_2. \tag{14}$$

Iterating the above process by using (14), we can obtain an $\succcurlyeq$-increasing sequence $x_0, x_1, x_2, \ldots$, such that

$$x_{n+1} \in Tx_n \text{ and } x_n \preccurlyeq x_{n+1}, \text{ for } n = 0, 1, 2, \ldots. \tag{15}$$

Since $D$ is a closed inductive subset of $X$, the sequential chain $\{x_n\} \subseteq D$ has an $\succcurlyeq$-upper bound $y_0 \in D$. Hence $\{x_n\}$ is $\succcurlyeq$-bounded with an $\succcurlyeq$-upper bound $y_0$ and an $\succcurlyeq$-lower bound $x_0$. Since the space $(X, \succcurlyeq, \|\cdot\|)$ is normal, so it is bounded (in $\|\cdot\|$). From the $\delta$-compactness of $T$, it implies that $\cup_{n\geq 0} Tx_n$ is relatively compact. From (15), we have

$$\{x_n\} \subseteq \cup_{n\geq 0} Tx_n.$$

Then $\{x_n\}$ has a convergent subsequence $\{x_{m(i)}\}$ and a point $x \in D$, such that

$$x_{m(i)} \to x, \text{ as } m(i) \to \infty. \tag{16}$$

Since $\{x_{m(i)}\}$ is also $\succcurlyeq$-increasing and it has limit $x$, by (16) and Lemma 2.3 [8], we get

$$x = \vee \{x_{m(i)}\}. \tag{17}$$

By the fact that $\{x_{m(i)}\}$ is a convergent subsequence of $\{x_n\}$, from (17), it implies

$$x = \vee \{x_n\}. \tag{18}$$

Since $\{x_{m(i)}\}$ is a subsequence of $\{x_n\}$ and $\{x_n\}$ is $\succcurlyeq$-increasing, then, for any $n \geq m(1)$, there is $i \geq 1$ such that $m(i) \leq n < m(i+1)$ satisfying

$$x_n \succcurlyeq x_{m(i)}.$$

It implies

$$\theta \preccurlyeq x - x_n \preccurlyeq x - x_{m(i)}.$$

Then, we have

$$\|x-x_n\| \leq \lambda\|x - x_{m(i)}\|, \tag{19}$$

where $\lambda \geq 1$ is the normality constant of $(X, \succcurlyeq, \|\cdot\|)$. Since $\{x_{m(i)}\}$ is a convergent sequence, combining (16) and (19) yields that $\{x_n\}$ is a convergent sequence in $D$ satisfying

$$x_n \to x, \text{ as } n \to \infty. \tag{20}$$

By (15) and (20) and from Lemma 3.5, $x \in Tx$. That is, $x$ is a fixed point of $T$. We obtain

$$\tilde{F}(T) \neq \emptyset.$$

From (17), we have $x_0 \preccurlyeq x$. It implies

$$\tilde{F}(T) \cap [x_0) \neq \emptyset.$$

We next prove that $\tilde{F}(T)$ is chain complete. Let $\{x_\alpha\}$ be an arbitrary ($\succcurlyeq$-increasing with respect to the index order) chain in $\tilde{F}(T)$. By the inductive property of $D$, without lose the generality, we assume that $\{x_\alpha\}$ has a lower $\succcurlyeq$-bound $u$ and an upper $\succcurlyeq$-bound $v$ in $D$. (If $D$ is bi-inductive, both of a lower $\succcurlyeq$-bound $u$ and an upper $\succcurlyeq$-bound $v$ of $\{x_\alpha\}$ exist. Otherwise, we can pick an arbitrary element $u = x_\beta \in \{x_\alpha\}$ and consider the sub-chain of $\{x_\alpha : x_\alpha \succcurlyeq x_\beta\}$, that has a lower $\succcurlyeq$-bound $x_\beta$ and, by the inductive property of $D$, an upper $\succcurlyeq$-bound $v$ of $\{x_\alpha\}$ exists). Hence $\{x_\alpha\} \subseteq [u, v]$, which is $\succcurlyeq$-bounded. Since $(X, \succcurlyeq, \|\cdot\|)$ is normal, then $[u, v]$ is bounded, so is $\{x_\alpha\}$. Since $T$ is $\delta$-compact, then $T(\{x_\alpha\})$ is relatively compact. From $x_\alpha \in Tx_\alpha$, for all $\alpha$, it implies that $\{x_\alpha\} \subseteq T(\{x_\alpha\})$ and $\{x_\alpha\}$ is relatively compact. So $\{x_\alpha\}$ is separable. In case if there is $x_\kappa \in \{x_\alpha\}$ such that $x_\alpha \preccurlyeq x_\kappa$, for all $\alpha$, then $x_\kappa = \vee\{x_\alpha\}$ and the chain-complete property of $\tilde{F}(T)$ is proved. Otherwise, from the separable property and the relative compactness of $\{x_\alpha\}$, we can select a sequence $\{x_m\} \subseteq \{x_\alpha\}$ such that

  (i) $\{x_m\}$ is $\succcurlyeq$-increasing;
  (ii) For every $x_\gamma \in \{x_\alpha\}$, there is $x_p \in \{x_m\}$ such that $x_\gamma \preccurlyeq x_p$;
  (iii) $\{x_m\}$ is a Cauchy sequence in $D$.

From (iii), there is $y \in D$, such that

$$x_m \to y, \text{ as } m \to \infty. \tag{21}$$

From (i) and (21), by Lemma 2.3 in [8], $y = \vee\{x_m\}$. By (ii), we get

$$y = \vee\{x_m\} = \vee\{x_\alpha\}. \tag{22}$$

Since $\{x_\alpha\} \subseteq \tilde{F}(T)$, it implies

$$x_m \in Tx_m, \text{ for } m = 1, 2, \ldots. \tag{23}$$

From (21) and (23) and by Lemma 3.5, we obtain $y \in Ty$. That is, $y \in \tilde{F}(T)$. By (22), it implies that $\tilde{F}(T)$ is chain-complete.

Part (b) can be similarly proved. Parts (a′) and (b′) are immediately consequences of part (a) and (b), respectively. □

**Corollary 3.7**. *Let* $(X, \|\cdot\|, \succcurlyeq)$ *be a normal partially ordered Banach space. Let* $u, v \in X$ *with* $u \prec v$. *Let* $T: [u, v] \to 2^{[u, v]} \setminus \{\emptyset\}$ *be an isotone and* $\delta$-*compact mapping with closed values. Then*

(b) $\mathcal{F}(T)$ *is a nonempty bi-chain-complete subset of* $[u, v]$;

(b′) *T has both* $\succcurlyeq$-*maximal and* $\succcurlyeq$-*minimal fixed points*.

*Proof*. From $u \prec v$, $[u, v]$ is an $\succcurlyeq$-interval, that is $\succcurlyeq$-closed and norm closed in $X$. Since $(X, \|\cdot\|, \succcurlyeq)$ is normal, then $[u, v]$ is bi-inductive and $\|\cdot\|$-bounded. Hence the values of $T$ are $\|\cdot\|$-bounded. It is clear to see that

$$u \preccurlyeq z, \text{ for all } z \in Tu.$$

So all conditions in Theorem 3.6 are satisfied and this corollary follows from Theorem 3.6 immediately. □

From Observation 3.4, we have the following consequence of Theorem 3.6.

**Corollary 3.8**. *Let* $(X, \|\cdot\|, \succcurlyeq)$ *be a normal partially ordered Banach space and D a closed inductive subset of X. Let* $F: D \to D$ *be an* $\succcurlyeq$-*increasing and compact mapping. Suppose that there is a point* $x_0 \in D$ *satisfying* $x_0 \preccurlyeq Fx_0$. *Then*

(a) $\mathcal{F}(F)$ *is a nonempty chain-complete subset of D*;

(b) *if D is bi-inductive, then* $\mathcal{F}(F)$ *is a nonempty bi-chain-complete subset of D*.

*Moreover*,

(a′) *F has an* $\succcurlyeq$-*maximal fixed point*;

(b′) *if D is bi-inductive, then F has both* $\succcurlyeq$-*maximal and* $\succcurlyeq$-*minimal fixed points*.

## 4. Fixed point theorems of single-valued mappings on regular partially ordered Banach spaces

Let $D$ be a closed subset of a Banach space $X$. A mapping $F: D \to X$ is said to be demi-continuous, whenever, for any given sequence $\{x_m\} \subseteq D$ and $x \in D$,

$$x_m \to x, \text{ strongly, as } m \to \infty \implies Fx_m \to Fx, \text{ weakly, as } m \to \infty.$$

$F: D \to X$ is said to be compact if $F$ is continuous and $F$ maps bounded sets to relatively compact sets.

**Theorem 4.1**. *Let* $(X, \|\cdot\|, \succcurlyeq)$ *be a regular partially ordered Banach space and D a closed inductive subset of X. Let* $F: D \to D$ *be an* $\succcurlyeq$-*increasing and demi-continuous mapping. Suppose that there is a point* $x_0 \in D$ *satisfying* $x_0 \preccurlyeq Fx_0$. *Then*

(a) $\mathcal{F}(F)$ *is a nonempty chain-complete subset of D*;

(b) *If D is bi-inductive, then $\bar{F}(F)$ is a nonempty bi-chain-complete subset of D.*

*Moreover,*

(a′) *F has $\succcurlyeq$-maximal fixed point;*

(b′) *if D is bi-inductive, then F has both $\succcurlyeq$-maximal and $\succcurlyeq$-minimal fixed points.*

*Proof.* For the given point $x_0 \in D$ in this theorem, denote $x_1 = Fx_0$. Then $x_0 \preccurlyeq x_1$. By the $\succcurlyeq$-increasing property of $F$, we can iteratively select an $\succcurlyeq$-increasing sequence $\{x_n\}$ such that

$$x_n \preccurlyeq x_{n+1} = Fx_n, \text{ for } n = 0, 1, 2, \ldots. \tag{24}$$

Since $D$ is a closed inductive subset of $X$, then the sequential chain $\{x_n\}$ has an upper $\succcurlyeq$-bound $y \in D$. From the regularity of $(X, \|\cdot\|, \succcurlyeq)$, $\{x_n\}$ contains a convergent subsequence $\{x_{m(i)}\}$ and a point $x \in D$, such that

$$x_{m(i)} \to x, \text{ as } m(i) \to \infty.$$

Notice that the regularity of partially ordered Banach spaces implies the normality. Similarly, to the proof of (18) and (20) in the proof of Theorem 3.6, we can show that

$$x = \vee\{x_n\}.$$

and

$$x_n \to x, \text{ as } n \to \infty. \tag{25}$$

Since $F$ is demi-continuous, by (24) and (25), we get

$$x_{n+1} = Fx_n \to Fx, \text{ weakly, as } n \to \infty. \tag{26}$$

(25) implies

$$x_{n+1} \to x, \text{ weakly, as } n \to \infty. \tag{27}$$

Combining (26) and (27) and by the uniqueness of weak limit of sequence, we get $Fx = x$. Hence $\bar{F}(T) \neq \emptyset$.

Next, we prove that $\bar{F}(T)$ is chain-complete. Take an arbitrary $\succcurlyeq$-increasing chain $\{y_\alpha\} \subseteq \bar{F}(T) \subseteq D$. Since regular partially ordered Banach spaces have chain-complete property, then $D$ is chain-complete. Then $\vee\{y_\alpha\}$ exists such that $\vee\{y_\alpha\} \in D$.

For any fixed $y_\beta \in \{y_\alpha\} \subseteq D$, by the regularity of $(X, \|\cdot\|, \succcurlyeq)$, we can show that $\{y_\alpha: y_\alpha \succcurlyeq y_\beta\}$ is relatively compact. So it is separable. Similarly, to the proof of (21) in the proof of Theorem 3.6, in case if there is $y_\kappa \in \{y_\alpha\}$ such that $y_\alpha \preccurlyeq y_\kappa$, for all $\alpha$, then $y_\kappa = \vee\{y_\alpha\}$ and the chain-complete property of $\bar{F}(T)$ is proved. Otherwise, from the separable property and the relatively compactness of $\{y_\alpha: y_\alpha \succcurlyeq y_\beta\}$, we can select a sequence $\{y_m\} \subseteq \{y_\alpha: y_\alpha \succcurlyeq y_\beta\}$ such that

(i) $\{y_m\}$ is $\succcurlyeq$-increasing;
(ii) for every $y_\gamma \in \{y_\alpha\}$, there is $y_p \in \{y_m\}$ such that $y_\gamma \preccurlyeq y_p$;
(iii) $\{y_m\}$ is a Cauchy sequence in $D$.

From (iii), there is $z \in D$, such that

$$y_m \to z, \text{ as } m \to \infty. \tag{28}$$

From Lemma 2.3 [8], we have

$$z = \vee\{y_n\}.$$

From properties (i) and (ii) of the selected sequence $\{y_m\}$, we have

$$z = \vee\{y_\alpha\}.$$

Since $Fy_\alpha = y_\alpha$, for all $\alpha$, from (28), we get

$$Fy_m \to Fz, \text{ weakly, as } m \to \infty.$$

That is,

$$y_m \to Fz, \text{ weakly, as } m \to \infty. \tag{29}$$

Combining (28) and (29) implies $Fz = z$. Hence

$$\vee\{y_\alpha\} = z \in \mathcal{F}(T).$$

It follows that $\mathcal{F}(T)$ is chain-complete. Part (b) can be similarly proved. Parts (a′) and (b′) are immediately consequences of part (a) and (b), respectively. □

## 5. Fixed point theorems of order-decreasing mappings on cones in partially ordered Banach space.

Let $(X, \succcurlyeq)$, $(U, \succcurlyeq^U)$ be posets and $F: X \to U$ a single-valued mapping. $F$ is said to be $\succcurlyeq$-decreasing whenever, for $x, y \in X$, $x \preccurlyeq y$ implies $F(x) \succcurlyeq^U F(y)$. An $\succcurlyeq$-increasing mapping $F: X \to U$ is said to be strictly $\succcurlyeq$-decreasing whenever $x \prec y$ implies $F(x) >^U F(y)$.

Let $(X, \|\cdot\|, \succcurlyeq)$ be a partially ordered Banach space, in which the partial order $\succcurlyeq$ is induced by a closed convex cone $K$ in $X$. Let $F: K \to K$ be a mapping. There are three conditions $H_1$, $H_2$, and $H_3$ listed in [5] for mapping $F$. They satisfy that $H_2 \Longrightarrow H_1$ and $H_3 \Longrightarrow H_1$. We only recall condition $(H_1)$ here:

$(H_1)$ $F: K \to K$ satisfies that $F^2\theta \succcurlyeq \epsilon F\theta$, for some $\epsilon \in (0, 1)$, and for any $\epsilon F\theta \preccurlyeq x \preccurlyeq F\theta$, $\epsilon \le t < 1$, there is $\eta = \eta(x, t)$ such that

$$F(tx) \preccurlyeq (t(1+\eta))^{-1} F(x).$$

We introduce the following conditions:

$(H^1)$ for $u, v \in K$, if $F(u) = v$ and $F(v) = u$, then $u = v$;
$(H^2)$ $\mathcal{F}(F^2) = \mathcal{F}(F)$.

**Lemma 5.1**. *Conditions $(H^1)$ and $(H^2)$ are equivalent*.

*Proof.* (H$^2$) $\Longrightarrow$(H$^1$). For $u, v \in K$, if $F(u) = v$ and $F(v) = u$, then $F^2(u) = u$ and $F^2(v) = v$. So $u, v \in \mathcal{F}(F^2)$. From (H$^2$), we have $u, v \in \mathcal{F}(F)$. It implies $u = F(u) = v$. Hence $u = v$.

(H$^1$) $\Longrightarrow$(H$^2$). Suppose $u \in \mathcal{F}(F^2)$. Let $F(u) = v$. Then $F(v) = F^2(u) = u$. From condition (H$^1$), it follows that $u = v$. It implies $u \in \mathcal{F}(F)$. The other inclusion $\mathcal{F}(F) \subseteq \mathcal{F}(F^2)$ is clear. $\square$

Theorem 3.2.1 in [5] is proved with condition (H$_1$). From Lemma 3.2.1 in [5], we see that (H$_1$)$\Longrightarrow$(H$^1$). In this section, by using the condition (H$^1$), we prove the following fixed point theorem, which can be considered as an extension of Theorem 3.2.1 in [5].

**Theorem 3.2.1** [5]. *Let $(X, \|\cdot\|, \succcurlyeq)$ be a partially ordered Banach space, in which the partial order is induced by a closed convex cone $K$ in $X$. Suppose that $F: K \to K$ is an $\succcurlyeq$-decreasing mapping satisfying condition* (H$_1$). *If one of the following two conditions satisfies*:

    (a) $(X, \|\cdot\|, \succcurlyeq)$ *is normal and $F$ is compact*;
    (b) $(X, \|\cdot\|, \succcurlyeq)$ *is regular and $F$ is demi-continuous*,

*then $F$ has a unique fixed point $x^*$ satisfying $x^* \preccurlyeq F\theta$. Furthermore, for every point $y_0 \in [\theta, F\theta]$,*

$$F^n y_0 \to x^*, \text{ as } n \to \infty. \tag{30}$$

**Theorem 5.2.** *Let $(X, \|\cdot\|, \succcurlyeq)$ be a partially ordered Banach space, in which the partial order $\succcurlyeq$ is induced by a closed convex cone $K$ in $X$. Suppose that $F: K \to K$ is an $\succcurlyeq$-decreasing mapping satisfying condition* (H$^1$). *If one of the following two conditions is satisfied*:

    (a)    $(X, \|\cdot\|, \succcurlyeq)$ *is normal and $F$ is compact*;
    (b)    $(X, \|\cdot\|, \succcurlyeq)$ *is regular and $F$ is demi-continuous*,

*then $F$ has a unique fixed point $x^*$ satisfying $x^* \preccurlyeq F\theta$. Furthermore, for every point $z \in K$,*

$$F^n z \to x^*, \text{ as } n \to \infty. \tag{30}$$

*Proof.* The proof is very similar to the proof of Theorem 3.2.1 in [5]. We first prove case (a). Let $\theta$ denote the origin of $X$. In case if $F\theta = \theta$, by the $\succcurlyeq$-decreasing property of $F$, it implies that $F(K) = \{\theta\}$ and $\mathcal{F}(F) = \{\theta\}$. This theorem is proved.

Next we suppose that $\theta \prec F\theta$ and we prove the existence of fixed point of $F$. Denote $x_0 = \theta$, $x_1 = F\theta$ and $x_n = F^n \theta$, for $n = 1, 2, \ldots$. Since $F$ is from $K$ to $K$ and $\succcurlyeq$ is induced by $K$, then $\theta = x_0 \preccurlyeq x_2$. From the $\succcurlyeq$-decreasing property of $F$, it implies

$$x_1 \succcurlyeq x_3, \ x_2 \preccurlyeq x_4, \ x_3 \succcurlyeq x_5, \ x_4 \preccurlyeq x_6, \ x_5 \succcurlyeq x_7, \ \ldots.$$

On the other hand, from $\theta \prec F\theta = x_1$, we have

$$x_1 \succcurlyeq x_2, \ x_2 \preccurlyeq x_3, \ x_3 \succcurlyeq x_4, \ x_4 \preccurlyeq x_5, \ x_5 \succcurlyeq x_6, \ \ldots.$$

Then we get an $\succcurlyeq$-increasing sequence $\{x_{2n}\}$ and an $\succcurlyeq$-decreasing sequence $\{x_{2n+1}\}$ satisfying

$$x_0 \preccurlyeq x_2 \preccurlyeq x_4 \preccurlyeq \ldots \preccurlyeq x_{2n} \preccurlyeq \quad \ldots \preccurlyeq \ldots \preccurlyeq x_{2n+1} \preccurlyeq \ldots \preccurlyeq x_5 \preccurlyeq x_3 \preccurlyeq x_1. \tag{31}$$

Since $\{x_{2n+1}\} \subseteq [\theta, F\theta]$ that is $\succcurlyeq$-bounded, by the normality of $(X, \|\cdot\|, \succcurlyeq)$, $\{x_{2n+1}\}$ is bounded. From the compactness of $F$, it implies that $\{x_{2n}\} = F(\{x_{2n+1}\})$ is relatively compact. Then $\{x_{2n}\}$ has a convergent subsequence. Similarly, to the proof of Theorem 3.6, we can show that $\{x_{2n}\}$ is convergent. So

$$\{x_{2n}\} \to u, \text{ as } n \to \infty, \text{ for some } u \in [\theta, F\theta]. \tag{32}$$

Since $\{x_{2n+1}\} = F(\{x_{2n}\})$, we can similarly prove that $\{x_{2n+1}\}$ is convergent and

$$\{x_{2n+1}\} \to v, \text{ as } n \to \infty, \text{ for some } v \in [\theta, F\theta]. \tag{33}$$

By the continuity of $F$, from (32) and (33), it implies

$$Fv = u \text{ and } Fu = v. \tag{34}$$

By Lemma 2.3 in [8] and (31), we get

$$u = \vee\{x_{2n}\} \preccurlyeq \wedge\{x_{2n+1}\} = v.$$

From the condition of this theorem that $F$ satisfies condition (H$^1$) and from (34), it implies $v = u$. Let $u = x^*$. By (34), it implies that $x^*$ is a fixed point of $F$.

From (32) and (33), it follows that $\{x_n\}$ is a convergent sequence with the limit point $u = x^*$. That is,

$$\{x_n\} \to x^*, \text{ as } n \to \infty. \tag{35}$$

To show the uniqueness of fixed point of $F$, we prove that, for every point $z \in K$, $\{F^n z\}$ is a convergent sequence in $[\theta, F\theta]$ and satisfies

$$F^n z \to x^*, \text{ as } n \to \infty. \tag{30}$$

To prove that, we denote $y_n = F^{n+1} z$, for $n = 0, 1, 2, \ldots$. From $\theta = x_0 \preccurlyeq z$ and the $\succcurlyeq$-decreasing property of $F$, we have $\theta = x_0 \preccurlyeq y_0 = Fz \preccurlyeq x_1 = Fx_0$. It implies $x_1 \succcurlyeq y_1 \succcurlyeq x_2$. Iterating this process, we obtain

$$x_{2n} \preccurlyeq y_{2n} \preccurlyeq x_{2n+1} \text{ and } x_{2n-1} \succcurlyeq y_{2n-1} \succcurlyeq x_{2n}, \text{ for } n = 0, 1, 2, \ldots.$$

It implies

$$\theta \preccurlyeq y_{2n} - x_{2n} \preccurlyeq x_{2n+1} - x_{2n} \text{ and } x_{2n-1} - x_{2n} \succcurlyeq y_{2n-1} - x_{2n} \succcurlyeq \theta, \text{ for } n = 0, 1, 2, \ldots. \tag{36}$$

Let $\lambda$ be the normality constant of $(X, \|\cdot\|, \succcurlyeq)$, (36) implies

$$\|y_{2n} - x_{2n}\| \leq \lambda\|x_{2n+1} - x_{2n}\| \text{ and } \|y_{2n-1} - x_{2n}\| \leq \lambda\|x_{2n-1} - x_{2n}\|, \text{ for } n = 0, 1, 2, \ldots. \tag{37}$$

Combining (35) and (37) implies (30) and moreover, $\mathcal{F}(F) \subseteq [\theta, F\theta]$. Then, for every $w \in \mathcal{F}(F)$, $\{F^n w\}$ is a constant sequence, in which all terms are $w$. By (30), it implies $w = x^*$.

Similarly, to the proofs of Theorem 4.1 and the case of part (a), we can prove this theorem for case (b). $\square$

By considering Remarks 5.0, Theorem 3.2.1 in [5] immediately follows from Theorem 5.2.

## 6. Applications to nonlinear integral equations

In this section, we apply Theorem 5.2 to prove the existence and uniqueness of solutions to some nonlinear integral equations, which were studied in [5].

**Example 6.1.** Let $\Psi: [0, 1] \to (0, 1]$ and $R$ be a real continuous function defined on $[0, 1] \times [0, 1]$. Consider the following nonlinear integral equation

$$1 = \Psi(x) + \Psi(x) \int_0^1 \frac{R(x, y)}{x^2 - y^2} \Psi(x) dy, \text{ for } 0 \le x \le 1. \tag{38}$$

Suppose that

(i) if $x \ge y$, then $R(x, y) \ge 0$; if $x < y$, then $R(x, y) \le 0$;
(ii) there are $v, M > 0$, and a nonnegative and bounded function $S$ defined on $[0, 1] \times [0, 1]$ with

$$\lim_{x, y \to 0^+} \frac{R(x, y)}{x + y} < \infty,$$

such that $\quad |R(x, y)| \le M|x-y|^v S(x, y)$, for all $x, y \in [0, 1]$.

Then, equation (38) has a unique solution $\Psi^*$ with $0 < \Psi^*(x) \le 1$, for $0 \le x \le 1$.

*Proof.* Let $X = C[0, 1]$, and $K = \{\psi \in C[0, 1]: \psi(x) \ge 0, \text{ for } 0 \le x \le 1\}$. Let $\succcurlyeq$ be the partial order on $X$ induced by $K$. Then $(X, \|\cdot\|, \succcurlyeq)$ is a normal partially ordered Banach space. Let

$$\psi(x) = \frac{1}{\Psi(x)} - 1.$$

Then $0 < \Psi(x) \le 1$ is equivalent to $\psi(x) \ge 0$. The equation (38) can be equivalently converted to

$$\psi(x) = \int_0^1 \frac{R(x, y)}{x^2 - y^2} \frac{1}{1 + \psi(x)} dy, \text{ for } 0 \le x \le 1. \tag{39}$$

Define $F: K \to K$ by

$$(F\psi)(x) = \int_0^1 \frac{R(x, y)}{x^2 - y^2} \frac{1}{1 + \psi(x)} dy, \text{ for } \psi \in K.$$

Then $F: K \to K$ is an $\succcurlyeq$-decreasing and completely continuous mapping. So it is compact. If there are $\psi, \varphi \in K$ such that

$$(F\psi)(x) = \int_0^1 \frac{R(x, y)}{x^2 - y^2} \frac{1}{1 + \psi(x)} dy = \varphi(x),$$

and

$$(F\varphi)(x) = \int_0^1 \frac{R(x, y)}{x^2 - y^2} \frac{1}{1+\varphi(x)} dy = \psi(x),$$

then

$$(1+ \varphi(x))\psi(x) = \int_0^1 \frac{R(x, y)}{x^2 - y^2} dy = (1+ \psi(x))\varphi(x).$$

It implies that $\psi(x) = \varphi(x)$. Hence $F$ satisfies condition (H$^1$). By Theorem 5.2, $F$ has a unique fixed point $\psi^* \in K$. So $\psi^*$ is the unique solution to equation (39). Then the function

$$\Psi^*(x) = \frac{1}{1+\psi^*(x)}$$

is the unique solution to equation (38). □